# MAJORIZATION FRAMEWORK FOR BALANCED LATTICE DESIGNS

By Aijun Zhang,[1] Kai-Tai Fang[2], Runze Li[3] and Agus Sudjianto

*University of Michigan, Hong Kong Baptist University, Pennsylvania State University and Bank of America*

This paper aims to generalize and unify classical criteria for comparisons of balanced lattice designs, including fractional factorial designs, supersaturated designs and uniform designs. We present a general majorization framework for assessing designs, which includes a stringent criterion of majorization via pairwise coincidences and flexible surrogates via convex functions. Classical orthogonality, aberration and uniformity criteria are unified by choosing combinatorial and exponential kernels. A construction method is also sketched out.

**1. Introduction.** We consider three types of balanced lattice designs including the fractional factorial design (FFD) [5], the supersaturated design (SSD) [1, 16] and the uniform design (UD) [10]. These have been widely used in agriculture, industry, scientific investigations and computer experiments, since a good design cannot only reduce experimental cost but also provide more efficient parameter estimation. Among many criteria for optimum factor assignment, minimum aberration [12, 17, 24, 26] considers the confounding situation between treatment effects under ANOVA decomposition; $E(s^2)$ [1] and $\text{Ave}(\chi^2)$ [27] measure two-factor orthogonality combinatorially; the discrepancy [14] considers the estimation of the overall mean from a multivariate quadrature perspective. The criteria are derived from

Received August 2003; revised January 2005.
[1]Supported by Hong Kong Baptist University, where Zhang completed his MPhil degree in 2004.
[2]Supported in part by the Hong Kong RGC Grants RGC/HKBU 2044/02P and RGC/HKBU 200804 and by Hong Kong Baptist University FRG Grant FRG/03-04/II-711.
[3]Supported by NSF Grants DMS-03-48869 and CCF-04-30349.

*AMS 2000 subject classification.* 62K15.

*Key words and phrases.* Admissible, discrepancy, fractional factorial design, majorization, minimum aberration, separable convex, supersaturated design, uniform design.

different principles, which may confuse users seeking to choose a suitable criterion for a specific experiment. Two natural questions are whether these individual criteria are connected and further, whether they can be unified into a single framework.

This paper aims to establish such a framework using majorization techniques. Majorization theory is appealing not only for its simplicity in concept, but its usefulness in many diverse fields; see [19] for a complete account. It was used as a tool in the early study of Kiefer's optimality criteria for optimal designs by considering eigenvalues of the Fisher information matrix; see [2, 22] and references therein. However, there has been little application to lattice designs, until the recent work of Cheng and Mukerjee [3] and Cheng, Steinberg and Sun [4] on estimation capacity, as well as that of Fang and Zhang [11] on projection aberration. In this paper we apply majorization theory to pairwise coincidences of experimental runs in order to study minimum aberration, discrepancy and some supersaturated design criteria. The lower bounds of these criteria will be provided in a unified way.

The paper is organized as follows. In Section 2 the majorization framework is proposed for balanced lattice designs, together with a two-stage investigation scheme through an illustrative example. Section 3 is devoted to unifying classical criteria surveyed above. In Section 4 an algorithm will be sketched out for constructing designs under the new framework. Technical proofs are given in the Appendix. Throughout this paper we use $|u|$ to denote the cardinality of a set $u$. The function $\binom{x}{j} = 0$ if $x < j$ and $\frac{1}{j!}x(x-1)\cdots(x-j+1)$ otherwise. The Kronecker delta $\delta(x,y) = 1$ if $x = y$ and is 0 otherwise.

**2. Majorization framework.** Consider experiments of $s$ factors each having $q$ levels. A *lattice design* with $n$ runs is a set of $n$ points chosen from the *lattice space* $\mathcal{L}(q^s)$, the $s$-fold tensor product of the integer set $\{0, 1, \ldots, q-1\}$. Each coordinate of $\mathcal{L}(q^s)$ corresponds to a factor. It is *balanced* (or *U-type*) when the $q$ levels appear equally often for each factor. The set of balanced designs is written as $\mathscr{U}(n, q^s)$. Either the *fractional factorial design* with resolution-$(t+1)$ or the *orthogonal array* $\text{OA}(n, s, q, t)$ lies in $\mathscr{U}(n, q^s)$, provided that, for any $t$ columns, all the possible level combinations appear equally often. The *uniform designs* are constructed from $\mathscr{U}(n, q^s)$. The *orthogonal designs* have strength $t \geq 2$ and they are *saturated* if $n = 1 + s(q-1)$; otherwise, the orthogonality is not attainable, as in *supersaturated designs*. For design selection, let $\mathscr{D}(n, q^s)$ denote the space of competing designs, which is restricted in this paper to be either $\mathscr{U}(n, q^s)$ or its subset.

Of our primary interest is the coincidence measurement between lattice points, which, together with its counterpart Hamming distance, plays an important role in the studies of designs and codes. For any $\mathbf{x}, \mathbf{w} \in$ $\mathcal{L}(q^s)$, the coincidence $\beta(\mathbf{x}, \mathbf{w}) := \sum_{j=1}^{s} \delta_{x_j, w_j}$ in terms of the Kronecker delta. It follows that $\beta(\mathbf{x}, \mathbf{x}) = s$ and $\beta(\mathbf{x}, \mathbf{w}) = \beta(\mathbf{w}, \mathbf{x})$. For a lattice design $\mathbf{X}(n, q^s)$, written as an $n \times s$ matrix with entries $x_{ij}$ from $\{0, 1, \ldots, q-1\}$, define its *pairwise coincidence* (PC) vector $\boldsymbol{\beta}(\mathbf{X}) := (\beta_1, \beta_2, \ldots, \beta_m)'$ by collecting $\beta(\mathbf{x}_i, \mathbf{x}_k)$ for $1 \leq i < k \leq n$ consecutively, where $m \equiv n(n-1)/2$ and $\beta(\mathbf{x}_i, \mathbf{x}_k) \equiv \beta_{n(i-1)+k-i(i+1)/2}$. We call two lattice designs PC-*different* if their PC-vectors cannot be exchanged by permutation. For isomorphic designs that are equivalent after reordering runs, permuting coordinates or switching levels, they hold the same increasing order statistic of PC-vector. The PC-different designs are nonisomorphic. The PC-sum $\sum_{r=1}^{m} \beta_r := \sum_{i<k} \sum_{j=1}^{s} \delta_{x_{ij}, x_{kj}}$ remains invariant in both isomorphic and nonisomorphic balanced designs, by observing that $1 + \sum_{k \neq i} \delta_{x_{ij}, x_{kj}} = n/q$ for any $i, j$.

LEMMA 1. *For any* $\mathbf{X} \in \mathscr{U}(n, q^s)$, *its PC-sum is* $\frac{ns}{2}(\frac{n}{q} - 1)$.

Let us now briefly review the majorization theory of Marshall and Olkin [19]. For a nonnegative vector $\mathbf{x} \in \mathbb{R}_+^m$, denote its increasing order statistic by $x_{[1]} \leq x_{[2]} \leq \cdots \leq x_{[m]}$. We say $\mathbf{x}$ is *majorized* by $\mathbf{y}$ and write $\mathbf{x} \preceq \mathbf{y}$ if

$$(2.1) \quad \sum_{r=1}^{k} x_{[r]} \geq \sum_{r=1}^{k} y_{[r]}, \qquad k = 1, 2, \ldots, m-1 \quad \text{and} \quad \sum_{r=1}^{m} x_r = \sum_{r=1}^{m} y_r.$$

If there exists at least one strict inequality $\sum_{r=1}^{k} x_{[r]} > \sum_{r=1}^{k} y_{[r]}$ for some $k$, we write $\mathbf{x} \prec \mathbf{y}$ strictly. A real-valued function $\Psi : \mathbb{R}_+^m \to \mathbb{R}$ is called *Schur-convex* if $\Psi(\mathbf{x}) \leq \Psi(\mathbf{y})$ for every pair $\mathbf{x}, \mathbf{y} \in \mathbb{R}_+^m$ with $\mathbf{x} \preceq \mathbf{y}$. Necessarily, $\Psi(\mathbf{x})$ is symmetric in its arguments, that is, invariant under permuting $x_1, \ldots, x_m$. We are mainly interested in the following *separable convex* class of Schur-convex functions:

$$\Psi(\mathbf{x}) = \sum_{r=1}^{m} \psi(x_r), \qquad \psi \text{ is convex on } \mathbb{R}_+,$$

as well as their monotonic mapping $g(\Psi(\mathbf{x}))$ for some $g$. Hardy, Littlewood and Pólya (HLP) [13] derived the following equivalent condition; or see page 108 of [19].

LEMMA 2 (HLP). *The inequality* $\Psi(\mathbf{x}) \leq \Psi(\mathbf{y})$ *holds for all separable convex functions* $\Psi : \mathbb{R}^m \to \mathbb{R}$ *if and only if* $\mathbf{x} \preceq \mathbf{y}$.

Consider the PC-mean of any balanced design $\mathbf{X}(n, q^s)$, which is a constant $\bar{\beta} = \frac{s(n-q)}{q(n-1)}$ by Lemma 1. For integer-valued $\boldsymbol{\beta}(\mathbf{X})$ with length $m$, let

$$\overline{\boldsymbol{\beta}} \equiv (\underbrace{\bar{\beta}, \ldots, \bar{\beta}}_{m})' \quad \text{and} \quad \widetilde{\boldsymbol{\beta}} \equiv (\underbrace{\theta, \ldots, \theta}_{m(1-f)}, \underbrace{\theta+1, \ldots, \theta+1}_{mf})',$$

where $\theta$ and $f$ are the integral part and fractional part of $\bar{\beta}$, respectively. It is clear that $\overline{\boldsymbol{\beta}} \preceq \widetilde{\boldsymbol{\beta}} \preceq \boldsymbol{\beta}(\mathbf{X})$, where $\widetilde{\boldsymbol{\beta}}$ reduces to $\overline{\boldsymbol{\beta}}$ when $f = 0$. By Lemma 2 we have a generalized version of Lemma 5.2.1 of Dey and Mukerjee [5].

LEMMA 3. *For integers $\beta_1, \ldots, \beta_m$ with mean $\bar{\beta}$, any separable convex function $\sum_{r=1}^{m} \psi(\beta_r)$ has a tight lower bound $m(1-f)\psi(\theta) + mf\psi(\theta+1)$, where $\theta$ and $f$ are the integral part and fractional part of $\bar{\beta}$, respectively.*

Based on decision theory and majorization theory, we propose a general framework for balanced designs via their PC measurements.

DEFINTION 1 (Majorization framework). In the space $\mathscr{D} \subseteq \mathscr{U}(n, q^s)$ of competing designs, we define the following:

1. $\mathbf{X}$ is inadmissible if there exists $\mathbf{X}'$ s.t. their PC-vectors satisfy $\boldsymbol{\beta}(\mathbf{X}') \prec \boldsymbol{\beta}(\mathbf{X})$;
2. $\mathbf{X}$ is majorant if $\boldsymbol{\beta}(\mathbf{X}) \preceq \boldsymbol{\beta}(\mathbf{X}')$ for all $\mathbf{X}' \in \mathscr{D}$;
3. $\mathbf{X}$ is Schur-$\psi$ optimal if $\Psi(\mathbf{X}; \psi) \leq \Psi(\mathbf{X}'; \psi)$ for all $\mathbf{X}' \in \mathscr{D}$, where

$$\text{Schur-}\psi \text{ criterion} \quad \Psi(\mathbf{X}; \psi) := \sum_{r=1}^{m} \psi(\beta_r(\mathbf{X}))$$

is determined by a specified convex kernel function $\psi : \mathbb{R}_+ \to \mathbb{R}$.

The three parts in Definition 1 can be divided hierarchically into two stages of investigation, namely, *stringent majorization check* and *flexible Schur-convex comparison*. At the first stage, for competing designs in $\mathscr{D}(n, q^s)$, compute their PC-vectors with elements sorted in increasing order. Compare the cumulative summations in the sense of the majorization ordering (2.1). By Definition 1.1, any inadmissible design should be prohibited for experimentation; by Definition 1.2, the majorant design(s) if it exists is the winner and absolutely recommended; otherwise, we need Definition 1.3 and go to the second stage for comparing admissible designs. The first stage is stringent since majorization requires strong conditions between PC-vectors. At the second stage, specify a convex kernel and compute the Schur-$\psi$ value for each admissible design. Since the above *Schur-$\psi$ criterion* is single-valued, all the designs are pairwise comparable and able to be rank-ordered. For different specific purposes, it is very flexible to predefine kernels, as long as they are convex functions. In the next section we shall discuss how to choose suitable kernels for investigating the orthogonality, aberration and uniformity properties of designs. Now let us illustrate the ideas with an example and some toy convex kernels.

TABLE 1
*Example of* 27*-run three-level design, transposed* $\mathbf{X}(27, 3^8)$. *Each row represents a factor, for which each level appears nine times. The* 70 *4-factor sub-designs are of interest*

| $U(27,3^8)$ | 1 | 2 | 3 | 4 | 5 | 6 | 7 | 8 | 9 | 10 | 11 | 12 | 13 | 14 | 15 | 16 | 17 | 18 | 19 | 20 | 21 | 22 | 23 | 24 | 25 | 26 | 27 |
|---|---|---|---|---|---|---|---|---|---|---|---|---|---|---|---|---|---|---|---|---|---|---|---|---|---|---|---|
| A | 1 | 2 | 0 | 2 | 0 | 0 | 0 | 1 | 2 | 1 | 0 | 2 | 2 | 1 | 0 | 2 | 0 | 2 | 0 | 1 | 1 | 1 | 2 | 1 | 1 | 0 | 2 |
| B | 1 | 2 | 0 | 0 | 0 | 2 | 2 | 2 | 1 | 0 | 2 | 1 | 2 | 0 | 1 | 0 | 1 | 1 | 1 | 1 | 2 | 2 | 2 | 1 | 0 | 0 | 0 |
| C | 2 | 2 | 1 | 1 | 1 | 1 | 0 | 1 | 0 | 0 | 2 | 1 | 0 | 0 | 0 | 2 | 2 | 1 | 0 | 1 | 2 | 1 | 0 | 2 | 0 | 2 | 2 |
| D | 1 | 1 | 2 | 0 | 1 | 2 | 1 | 0 | 2 | 1 | 0 | 2 | 2 | 2 | 0 | 1 | 0 | 1 | 1 | 0 | 2 | 1 | 0 | 2 | 0 | 2 | 0 |
| E | 1 | 2 | 1 | 1 | 2 | 2 | 0 | 0 | 2 | 0 | 2 | 0 | 1 | 2 | 2 | 2 | 0 | 0 | 1 | 2 | 1 | 1 | 1 | 0 | 0 | 0 | 1 |
| F | 0 | 0 | 0 | 1 | 0 | 2 | 1 | 0 | 1 | 2 | 1 | 1 | 0 | 1 | 0 | 2 | 2 | 2 | 2 | 2 | 2 | 1 | 2 | 0 | 0 | 1 | 1 |
| G | 1 | 2 | 2 | 2 | 0 | 1 | 2 | 0 | 0 | 0 | 0 | 0 | 1 | 2 | 1 | 1 | 2 | 1 | 0 | 2 | 0 | 1 | 2 | 2 | 1 | 1 | 0 |
| H | 2 | 0 | 2 | 0 | 1 | 1 | 2 | 0 | 1 | 2 | 2 | 2 | 2 | 0 | 0 | 2 | 1 | 0 | 0 | 2 | 0 | 1 | 1 | 1 | 1 | 0 | 1 |

EXAMPLE 1. For the purpose of illustration, consider the following scenario: an experiment of 27 runs with 8 factors each having 3 levels. The experimenter chooses a uniform design tabulated in Table 1, which was obtained by Fang, Ma and Winker [8]. The experimenter has some prior knowledge: among the 8 factors, 4 factors may have potential impact on the output, while the other 4 have little impact; he is interested in including all of them in the study. To incorporate such prior information into the design of the experiment, he wants to choose a sub-design consisting of 27 runs with 4 factors for the 4 potential factors. This leads us to study the following problem: how to choose the 4-factor sub-design from Table 1?

There are in total $\binom{8}{4} = 70$ choices of sub-designs from this table, which are all balanced and form the design space $\mathscr{D}(27, 3^4)$ of this study. For demonstration, 4 sub-designs from the design space $\mathscr{D}$ are chosen through $\mathbf{X}_1 = \{A, C, G, H\}$, $\mathbf{X}_2 = \{B, C, G, H\}$, $\mathbf{X}_3 = \{A, B, D, F\}$ and $\mathbf{X}_4 = \{A, D, E, F\}$ labels of factors. Let us make a two-stage investigation of $\mathscr{D}$ under the majorization framework:

• Stage 1: stringent majorization check.

For each sub-design, its PC-vector has length 351 and sum 972. By majorization ordering, not all 70 sub-designs can be compared. There exists

TABLE 2
*Numerical results of Schur-convex comparison for* $\mathbf{X}_1$ *to* $\mathbf{X}_4$ *from* $\mathscr{D}(27, 3^4)$. *The lower bounds in the last column are derived from Theorem* 1 (*where* $\sum$ *denotes* $\sum_{r=1}^{m}$)

| Convex kernel | $\Psi(\mathbf{X}; \psi)$ | $\mathbf{X}_1$ | $\mathbf{X}_2$ | $\mathbf{X}_3$ | $\mathbf{X}_4$ | Lower bound |
|---|---|---|---|---|---|---|
| Variance | $\frac{1}{m}\sum(\beta_r - \bar{\beta})^2$ | 0.6391 | 0.6391 | 0.6732 | 0.6789 | 0.1775 |
| Power | $\sum \beta_r^\pi$ | 1658.7 | 1724.5 | 1765.5 | 1790.4 | 984.8 |
| Exponential | $\sum(\frac{1+\sqrt{5}}{2})^{\beta_r}$ | 683.4 | 685.6 | 687.9 | 688.5 | 648.9 |

no majorant design in the given $\mathscr{D}$. For $\mathbf{X}_1$ to $\mathbf{X}_4$, we have

$$\boldsymbol{\beta}(\mathbf{X}_1) \prec \boldsymbol{\beta}(\mathbf{X}_3) \prec \boldsymbol{\beta}(\mathbf{X}_4), \qquad \boldsymbol{\beta}(\mathbf{X}_2) \prec \boldsymbol{\beta}(\mathbf{X}_3) \prec \boldsymbol{\beta}(\mathbf{X}_4),$$

where both $\mathbf{X}_3$ and $\mathbf{X}_4$ are inadmissible (even though $\mathbf{X}_3$ overwhelms $\mathbf{X}_4$), but the admissible $\mathbf{X}_1$ and $\mathbf{X}_2$ are not distinguishable at this stage.
- Stage 2: Schur-convex comparison.

Let us choose three toy kernels for comparisons, namely, a variance kernel $\psi_1(x) = m^{-1}(x - \bar{x})^2$, a $\pi$th-power kernel $\psi_2(x) = x^\pi$ and an exponential kernel based on the golden ratio $\psi_3(x) = (\frac{1+\sqrt{5}}{2})^x$. Numerical results of 4 sub-designs are shown in Table 2. Their Schur-$\psi$ values are rank-ordered as

$$\Psi(\mathbf{X}_1; \psi_j) \leq \Psi(\mathbf{X}_2; \psi_j) < \Psi(\mathbf{X}_3; \psi_j) < \Psi(\mathbf{X}_4; \psi_j) \qquad \text{for } j = 1, 2, 3,$$

where the equality holds only for $\psi_1$. It is shown that inadmissible $\mathbf{X}_3, \mathbf{X}_4$ always have large Schur-$\psi$ values no matter what convex kernel is used.

Under the classical criteria, both $\mathbf{X}_1$ and $\mathbf{X}_2$ are orthogonal designs of resolution 3; their (generalized) word-length patterns are given by

$$\mathrm{GWP}(\mathbf{X}_1) = (0, 0, 10/9, 8/9), \qquad \mathrm{GWP}(\mathbf{X}_2) = (0, 0, 46/27, 20/27)$$

and the wrap-around $L_2$-discrepancy values are given by $\mathrm{WL}_2(\mathbf{X}_1) = 0.4242$ and $\mathrm{WL}_2(\mathbf{X}_2) = 0.4245$. We find that in the complete pool $\mathscr{D}(27, 3^4)$ of 70 competing designs, $\mathbf{X}_1$ is not only an FFD with minimum aberration but also a UD with minimum $\mathrm{WL}_2$-discrepancy.

The above example demonstrates both stringency and flexibility of the majorization framework for assessing designs. The kernel selection problem at stage 2 is discussed in some detail by Zhang [28], who also explains why $\mathbf{X}_1$ and $\mathbf{X}_2$ are not distinguishable under the variance kernel. Formally, we have the following main theorem to characterize the necessary and sufficient conditions between majorant designs and Schur-$\psi$ optimum designs, according to Lemma 2. We also employ Lemma 3 to derive the lower bounds for specific Schur-$\psi$ criteria.

THEOREM 1. *A balanced lattice design is majorant if and only if it is Schur-$\psi$ optimum w.r.t. every convex kernel. For any well-defined Schur-$\psi$ criterion, it has a lower bound $m(1-f)\psi(\theta) + mf\psi(\theta+1)$.*

The lower bound is presented for general PC-mean $\bar{\beta}$, either integer-valued or not. Obviously, if $\frac{s(n-q)}{q(n-1)}$ is a positive integer, $f = 0$ and the lower bound reduces to $m\psi(\bar{\beta})$. This bound is attainable if there exists an *equidistant design* $\overline{\mathbf{X}}$ in $\mathscr{U}(n, q^s)$ such that all the Hamming distances between distinct runs are identical, that is, $\mathrm{PC}(\overline{\mathbf{X}}) = \overline{\boldsymbol{\beta}}$. Equidistant designs are a typical type of majorant design, and examples are two-level SSDs constructed by the half-fraction Hadamard method [16], multi-level SSDs constructed from resolvable balanced incomplete block designs [6] and saturated $OA(n, s, q, 2)$ designs whose $\beta(\mathbf{x}_i, \mathbf{x}_k) \equiv s - n/q$ for any $i \neq k$ [20]. If $\frac{s(n-q)}{q(n-1)}$ returns a noninteger, the lower bound by Theorem 1 can be achieved by *weak equidistant designs* with $\widetilde{\boldsymbol{\beta}}$ whose elements differ at most by 1. Examples can be obtained by either adding a balanced factor to or removing a factor from saturated designs. Note that the bound is tight in some cases but not generally tight under all parameter $(n, s, q)$ settings.

**3. Unification of classical criteria.** The design criteria for FFD, SSD and UD are discussed in this section. The majorization framework and, in particular, the flexible Schur-$\psi$ criteria based on combinatorial and exponential kernels are used to unify the criteria of minimum aberration and discrepancy. By Theorem 1, their lower bounds are generated automatically. Throughout this section $\theta$ is the integral part of $\bar{\beta}$ [in particular, $\bar{\beta} = \frac{s(n-q)}{q(n-1)}$ for $\mathbf{X}(n, q^s)$] and $f = \bar{\beta} - \theta$.

3.1. *Fractional factorial designs.* FFD is an important experimental strategy and usually measured by the minimum aberration criterion originally proposed by Fries and Hunter [12] for regular designs. We rely on its generalization by Xu and Wu [26] for both two- and multi-level, both regular and nonregular designs. Based on the ANOVA decomposition model, define for $\mathbf{X}(n, q^s)$

$$(3.1) \qquad A_j(\mathbf{X}) := \frac{1}{n^2} \|\mathbf{G}_j\|_F^2 \equiv \frac{1}{n^2} \text{trace}(\mathbf{G}_j^H \mathbf{G}_j), \qquad j = 1, \ldots, s,$$

where $\mathbf{G}_j$ is the matrix consisting of all $j$-factor contrast coefficients ($\|\cdot\|_F$: Frobenius norm; $^H$: conjugate transpose). The (generalized) word-length pattern (GWP) is defined by $(A_1, \ldots, A_s)$, in which $A_1 \equiv 0$ for balanced designs. For two such patterns $\mathbf{x}, \mathbf{y} \in \mathbb{R}_+^s$, define a partial ordering $\models$ as follows. We write $\mathbf{x} \vdash \mathbf{y}$ if the first nonzero element of $\mathbf{x} - \mathbf{y}$ is negative, and write $\mathbf{x} \models \mathbf{y}$ if $\mathbf{x} \vdash \mathbf{y}$ or $\mathbf{x} = \mathbf{y}$. An FFD has minimum aberration if its GWP achieves the minimum under $\models$. Ma and Fang [17] and Xu and Wu [26] connected word-length pattern with MacWilliams' transform of distance distribution in coding theory,

$$(3.2) \qquad A_j(\mathbf{X}) = \frac{1}{n} \sum_{l=0}^{s} E_l(\mathbf{X}) P_j(l; s, q), \qquad j = 1, \ldots, s,$$

where $E_l(\mathbf{X}) = n^{-1} |\{(\mathbf{x}_i, \mathbf{x}_k) : \beta(\mathbf{x}_i, \mathbf{x}_k) = s - l, \ i, k = 1, \ldots, n\}|$ for $l = 0, \ldots, s$ and

$$P_j(x; s, q) = \sum_{w=0}^{j} (-1)^w (q-1)^{j-w} \binom{x}{w} \binom{s-x}{j-w}$$

are *Krawtchouk polynomials* ([18], Section 5.7). Clearly, $A_j(\mathbf{X})$ can be expressed as $\frac{2}{n^2}\sum_{r=1}^{m} P_j(s-\beta_r;s,q) + \frac{(q-1)^j}{n}\binom{s}{j}$. To unify the minimum aberration through Schur-$\psi$ criterion, a direct idea is to use Krawtchouk polynomials, $P_2(s-x;s,q)$ to $P_s(s-x;s,q)$. However, the function $P_j(s-x;s,q)$ is not generally convex except for $j=2$, which implies that it is trivial to unify $A_2$ and find its lower bound, but nontrivial for higher-order $A_j$'s.

Let us make an indirect approach by a series of combinatorial functions. For $\mathbf{X}(n,q^s)$, define the *Schur-combinatorial criterion* of affine form,

$$(3.3) \quad \Psi_C(\mathbf{X};j) := 2\sum_{r=1}^{m}\binom{\beta_r}{j} - \binom{s}{j}\left(\frac{n^2}{q^j} - n\right), \qquad j=1,\ldots,s,$$

which are all separable convex on $\mathbb{R}_+^m$. The criterion can be interpreted statistically as follows. Consider $\Psi_C(\mathbf{X};s)$ first. Randomly on $\mathcal{L}(q^s)$, each $\mathbf{x}$ has uniform probability $n/q^s$ of entering the $n$-point design. Let $N_\mathbf{x}$ be its true occurrences in $\mathbf{X}(n,q^s)$. The value $(N_\mathbf{x} - n/q^s)^2$ measures the variation of the design centering $\mathbf{x}$. For $q^s$ different points, the total variation $\sum_{\mathbf{x}\in\mathcal{L}(q^s)}(N_\mathbf{x} - n/q^s)^2$ is therefore a measure of uniform covering, which equals $\Psi_C(\mathbf{X};s)$. Formally, we have:

THEOREM 2. *For $\mathbf{X}(n,q^s)$ and $\mathcal{S} = \{1,\ldots,s\}$, the Schur-combinatorial criterion*

$$(3.4) \quad \Psi_C(\mathbf{X};j) = \sum_{u\subseteq\mathcal{S},|u|=j}\sum_{\mathbf{x}\in\mathcal{L}(q^j)}\left(N_\mathbf{x}^{(u)} - \frac{n}{q^j}\right)^2 \qquad \text{for } j=1,\ldots,s,$$

*where $N_\mathbf{x}^{(u)}$ counts the runs whose $u$-coordinates take level-combination $\mathbf{x}$. Further, the design $\mathbf{X}$ has orthogonal strength $t$ if and only if $\Psi_C(\mathbf{X};j)=0$ for $j=1,\ldots,t$.*

Projection properties are taken into account in Theorem 2. Let $\mathbf{X}_u$ denote the $u$-coordinate sub-design. Thus, $\Psi_C(\mathbf{X};j) = \sum_{|u|=j}\Psi_C(\mathbf{X}_u;j)$ sums up the measurements at all $j$-dimensional sub-spaces. Besides the geometrical meaning, $\Psi_C(\mathbf{X};j)$ measures the variation from the $j$-factor orthogonal strength. From (3.4), $\Psi_C(\mathbf{X};j) \geq 0$, where equality holds if $N_\mathbf{x}^{(u)} \equiv n/q^j$, which occurs if and only if $\mathbf{X}$ is an orthogonal array of strength $j$.

To use a relatively simple notation, define for $\mathbf{X}(n,q^s)$ the root-mean-squared *deviation* criterion,

$$(3.5) \quad B_s(\mathbf{X}) := \sqrt{\frac{1}{q^s}\sum_{\mathbf{x}\in\mathcal{L}(q^s)}(N_\mathbf{x} - n/q^s)^2} \equiv \sqrt{\frac{1}{q^s}\Psi_C(\mathbf{X};s)},$$

as well as $B_j(\mathbf{X}) := \sqrt{\frac{1}{q^j}\Psi_C(\mathbf{X};j)}$ for $j<s$. Let us call $(B_1,\ldots,B_s)$ a *deviation pattern*, which reduces to the projection $V$-criterion for two-level designs [23]. Note that $B_t = 0$ implies that $B_j = 0$ for $j < t$. Analogous to the word-length pattern, $A_1 = B_1 = 0$ for balanced lattice designs, and $A_t = D_t = 0$, $A_{t+1} > 0$, $B_{t+1} > 0$ for resolution-$(t+1)$ orthogonal designs.

THEOREM 3. *For design* $\mathbf{X} \in \mathscr{U}(n, q^s)$, *the deviation pattern and word-length pattern are linearly related by*

$$(3.6a) \qquad B_j^2(\mathbf{X}) = \frac{n^2}{q^{2j}} \sum_{k=1}^{j} \binom{s-k}{j-k} A_k(\mathbf{X}) \qquad \text{for } j = 1, \ldots, s.$$

*Their benchmarks are given by* $(0, A_2^*, \ldots, A_s^*)$ *and* $(0, B_2^*, \ldots, B_s^*)$, *in which*

$$(3.6b) \qquad \begin{aligned} A_j^* &= \left(1 - \frac{1}{n}\right)((1-f)P_j(s-\theta; s, q) + fP_j(s-\theta-1; s, q)) \\ &\quad + \frac{(q-1)^j}{n}\binom{s}{j}, \end{aligned}$$

$$(3.6c) \quad B_j^* = \sqrt{\frac{n(n-1)}{q^j}\left(\binom{\theta}{j} + f\binom{\theta}{j-1}\right) - \binom{s}{j}\left(\frac{n^2}{q^{2j}} - \frac{n}{q^j}\right)},$$

*for* $j = 2, \ldots, s$, *in the sense that* $(0, A_2^*, \ldots, A_s^*) \models (0, A_2(\mathbf{X}), \ldots, A_s(\mathbf{X}))$ *while* $B_j^* \leq B_j(\mathbf{X})$ *for all* $j$ *simultaneously*.

Checking the simplest case for balanced designs with integer-valued PC-mean $\bar{\beta}$, we get $A_2(\mathbf{X}(n, q^s)) \geq \frac{s(q-1)(qs-s-n+1)}{2(n-1)}$, which is consistent with Fang, Ge, Liu and Qin [6] and Xu [25] for investigating supersaturated designs.

3.2. *Supersaturated designs.* In the recent decade, SSDs, in most cases 2-level factorials, have drawn much attention in screening experimentation due to their economic run size. Nonorthogonality criteria like $E(s^2)$ and $\text{Ave}(\chi^2)$ are used to evaluate/construct SSDs, as their orthogonal property is violated. For $\mathbf{X}(n, 2^s)$, let $\mathbf{x}_{(j)}$ with $(-1, 1)$ entries represent the $j$th factor. Booth and Cox [1] originally defined $E(s^2)$ by the mean inner-product $\frac{2}{s(s-1)} \sum_{1 \leq j < l \leq s} x_{(j)}^T \mathbf{x}_{(l)}$. Let $N_{\tau_1, \tau_2}^{(j,l)}$ be the number of runs whose $(j, l)$ factors take level-combination $(\tau_1, \tau_2)$. Then we observe that $\mathbf{x}_{(j)}^T \mathbf{x}_{(l)} = 4 \sum_{\tau_1, \tau_2 = 1}^{2} (N_{\tau_1, \tau_2}^{(j,l)} - n/4)^2$. For multi-level $\mathbf{X}(n, q^s)$, define

$$\text{Ave}(\chi^2) := \frac{2}{s(s-1)} \sum_{1 \leq j < l \leq s} \sum_{\tau_1, \tau_2 = 1}^{q} (N_{\tau_1, \tau_2}^{(j,l)} - n/q^2)^2,$$

which reduces to $E(s^2)$ when $q = 2$ (after being multiplied by 4) and reduces to Yamada and Lin's [27] $\text{Ave}(\chi^2)$ when $q = 3$ (after being multiplied by $9/n$).

By Theorem 2, we find that there is a natural link between $\text{Ave}(\chi^2)$ and $\Psi_C(\mathbf{X}; 2)$ based on the combinatorial kernel, as well as the deviation measure $B_2$ on the 2D sub-space. For simplicity, we give a unification scheme through the quadratic kernel $\psi(\beta) = \beta^2$ and the associated Schur-$\psi$ criterion $\Psi(\mathbf{X}; \beta^2) := \sum_{r=1}^{m} \beta_r^2$, which is equivalent to $\Psi_C(\mathbf{X}; 2)$ for balanced designs.

THEOREM 4. *The nonorthogonality criterion* $\text{Ave}(\chi^2)$ *for SSD* $\mathbf{X} \in \mathscr{U}(n, q^s)$ *satisfies*

$$\text{Ave}(\chi^2) = \frac{2}{s(s-1)} \Psi(\mathbf{X}; \beta^2) + a \geq \frac{n(n-1)}{s(s-1)} (\theta^2 + 2\theta f + f) + a,$$

*where the constant* $a = \frac{q^2 n s + n^2(1-s-q)}{q^2(s-1)}$.

The lower bound follows directly Theorem 1. When the PC-mean $\bar{\beta}$ is an integer,

$$E(s^2) \geq \frac{n^2(s-n+1)}{(s-1)(n-1)}, \qquad \text{Ave}(\chi^2) \geq \frac{n^2(q-1)((q-1)s - n - 1)}{q^2(s-1)(n-1)},$$

where the lower bounds can be attained by optimum SSDs constructed from partial saturated designs, resolvable BIBDs or an algorithmic approach [6, 16, 21, 27].

3.3. *Uniform designs.* UD is of space filling type and becomes more and more important for computer experiments. For $n$ design points scattered into the lattice space, Fang and Wang [10] suggested using the star discrepancy as the uniformity measure, which corresponds to the famous *Kolmogorov–Smirnov statistic* for goodness-of-fit testing between $F_n(\mathbf{x})$, the empirical distribution of the design, and $F_*(\mathbf{x})$, the uniform distribution. A discrepancy defined in quasi-Monte Carlo methods can be viewed as a norm $\|F_n(\mathbf{x}) - F_*(\mathbf{x})\|$ of some reproducing kernel Hilbert space [14]. For assessing qualitative factorial assignment, Hickernell and Liu [15] proposed the discrete discrepancy, which is a special case (when $\mu = 0$) of the *categorical discrepancy* below.

DEFINTION 2 (Categorical discrepancy). On the lattice space $\mathcal{L}(q^s)$ and the set $\mathcal{S} = \{1, \ldots, s\}$, for each nonempty $u \subseteq \mathcal{S}$, any $\mathbf{x}, \mathbf{w} \in \mathcal{L}(q^s)$ and any design $\mathbf{X}$ with points $\mathbf{x}_1, \ldots, \mathbf{x}_n \in \mathcal{L}(q^s)$, define the categorical type of hat reproducing kernel function and hat discrepancy

$$\widehat{\mathbb{K}}_u(\mathbf{x}, \mathbf{w}) = \prod_{j \in u} (b + (a-b)\delta_{x_j, w_j}),$$

$$D_u(\mathbf{X}; \widehat{\mathbb{K}}_u) = \left( -\mu^{|u|} + \frac{1}{n^2} \sum_{i,k=1}^{n} \widehat{\mathbb{K}}_u(\mathbf{x}_i, \mathbf{x}_k) \right)^{1/2},$$

where $a$ is a given positive constant $(a < q-1)$, $b$ is chosen from $[-\frac{a}{q-1}, a)$ and $\mu = \frac{1}{q}(a + (q-1)b)$. Define the categorical discrepancy pattern $(D_1, \ldots, D_s)$ and the categorical discrepancy $D(\mathbf{X}; a \vee b)$, respectively, by

$$D_j(\mathbf{X}; a \vee b) = \sqrt{\sum_{|u|=j} D_u^2(\mathbf{X}; \widehat{\mathbb{K}}_u)}$$

and

$$D(\mathbf{X}; a \vee b) = \sqrt{\sum_{j=1}^{s} D_j^2(\mathbf{X}; a \vee b)},$$

where $a \vee b$ denotes categorical assignments to the hat kernel $\widehat{\mathbb{K}}_u$.

The parameter constraints $a > 0$ and $-\frac{a}{q-1} \leq b < a$ are set to ensure that the bivariate $\widehat{\mathbb{K}}_u$ is nonnegative definite. For $j = 1, \ldots, s$, $D_j(\mathbf{X}; a \vee b)$ sums up the hat discrepancies of all possible $j$-factor projection designs. Hickernell and Liu [15] showed that when the parameters satisfy $a + (q-1)b = 0$, the categorical discrepancy pattern under partial ordering $\models$ is equivalent to the minimum aberration criterion. For $\mathbf{X}(n, q^s)$ under our majorization framework, by using the kernel $\psi(\beta) = \rho^\beta$ with base $\rho > 1$, we can define the *Schur-exponential criterion* $\Psi_\mathrm{E}(\mathbf{X}; \rho) := \sum_{r=1}^{m} \rho^{\beta_r}$.

THEOREM 5. *For any lattice design* $\mathbf{X} \in \mathscr{U}(n, q^s)$, *the squared categorical discrepancy is equivalent to the Schur-exponential criterion,*

$$(3.6\mathrm{g}) \qquad D^2(\mathbf{X}; a \vee b) = \frac{2\Psi_\mathrm{E}(\mathbf{X}; \rho)}{n^2} + \frac{(1+a)^s}{n} - (1+\mu)^s,$$

*where the exponential base* $\rho = (1+a)/(1+b)$. *It has lower bound*

$$D^2(\mathbf{X}; a \vee b) \geq \frac{1}{n}((n-1)(1 - f + \rho f)\rho^\theta + (1+a)^s) - (1+\mu)^s.$$

The *centered $L_2$-discrepancy* (CL$_2$) and *wrap-around $L_2$-discrepancy* (WL$_2$) are popular uniformity measures for quantitative experiments; see [14] and [8] for the details. For modest-level designs, CL$_2$- and WL$_2$-discrepancies have similar properties to categorical discrepancy, that is, the reproducing kernel values between distinct runs are determined by coincidence measurement. They correspond to the Schur-exponential criteria under different bases.

COROLLARY 1. *For* $\mathbf{X} \in \mathscr{U}(n, q^s)$, *the Schur-exponential criterion can cover*

$$q = 2 : \mathrm{CL}_2^2(\mathbf{X}) - a_1 = \frac{2\Psi_\mathrm{E}(\mathbf{X}; 1.25)}{n^2} \geq \frac{(n-1)(4+f)}{4n}\left(\frac{5}{4}\right)^\theta,$$

$$q = 2 : \mathrm{WL}_2^2(\mathbf{X}) - a_2 = \frac{2\Psi_\mathrm{E}(\mathbf{X}; 1.2)}{n^2} \left(\frac{5}{4}\right)^s \geq \frac{(n-1)(5+f)}{5n} \left(\frac{5}{4}\right)^s \left(\frac{6}{5}\right)^\theta,$$

$$q = 3 : \mathrm{WL}_2^2(\mathbf{X}) - a_2 = \frac{2\Psi_\mathrm{E}(\mathbf{X}; 27/23)}{n^2} \left(\frac{23}{18}\right)^s$$

$$\geq \frac{(n-1)(23+4f)}{23n} \left(\frac{23}{18}\right)^s \left(\frac{27}{23}\right)^\theta,$$

where $a_1 = \frac{1}{n}(\frac{5}{4})^s + (\frac{13}{12})^s - 2(\frac{35}{32})^s$ and $a_2 = \frac{1}{n}(\frac{3}{2})^s - (\frac{4}{3})^s$.

The lower bounds derived above for $\mathrm{CL}_2$-discrepancy and $\mathrm{WL}_2$-discrepancy are tighter than [9] and [6].

**4. Algorithmic construction.** In our framework, majorization on pairwise coincidences is conceptually simple. From a geometric point of view, it enforces pairwise coincidences spread as equally as possible, which is universally applicable to various criteria discussed above. By the majorization idea, an optimization approach with heuristic searches is in development. Although this paper addresses mainly theoretical aspects of design criteria, we briefly sketch out our algorithmic construction method, in particular, its basic operation of Robin Hood swap. The swapping algorithm aims to take one unit from the $\beta$-large pair of coincident runs and give it to the $\beta$-small pair, analogous to the legend of Robin Hood.

ALGORITHM 1 (Robin Hood swap pseudo-code). Given a convex kernel function $\psi : \mathbb{R}_+ \to \mathbb{R}$ and a balanced lattice design $\mathbf{X} \in \mathscr{U}(n, q^s)$:

Step 1: Compute the coincidence matrix $\mathbf{M}$ and find its maximal entry (entries); for each such pairwise coincidence $\beta_{ik}$ of runs $(\mathbf{x}_i, \mathbf{x}_k)$, do steps 2 and 3.

Step 2: Find the run(s) that has minimal coincidence from $\mathbf{x}_i$. For each such run $\mathbf{x}_t$, find the coordinates $\mathscr{C}$ such that $x_{ij} = x_{kj}$, while $x_{tj} \neq x_{ij}$ for $j \in \mathscr{C}$.

Step 3: For each coordinate $j \in \mathscr{C}$, find $\mathscr{R}_i, \mathscr{R}_t$ such that $x_{wj} = x_{ij}, \forall w \in \mathscr{R}_i$ and $x_{wj} = x_{tj}, \forall w \in \mathscr{R}_t$, respectively ($i \notin \mathscr{R}_i, t \notin \mathscr{R}_t$). Compute the delta,

$$\Delta_j = \sum_{w \in \mathscr{R}_i} (\psi(\beta_{iw} - 1) + \psi(\beta_{tw} + 1))$$

$$+ \sum_{w \in \mathscr{R}_t} (\psi(\beta_{tw} - 1) + \psi(\beta_{iw} + 1))$$

$$- \sum_{w \in \mathscr{R}_i \cup \mathscr{R}_t} (\psi(\beta_{iw}) + \psi(\beta_{tw})).$$

Find the local minimum $\Delta_*^{(i,t)}$ for $j \in \mathscr{C}$. Record $\{i, t, j_*; \Delta_*^{(i,t)}\}$ if $\Delta_*^{(i,t)} < 0$.

Step 4: Find from the record the global minimum $\Delta_*$ and output $\{i_*, t_*, j_*\}$.

The algorithm works on a specific kernel function. For example, consider the quadratic kernel $\psi(x) = x^2$ and the randomly generated balanced design $\mathbf{X} \in \mathscr{U}(8, 2^6)$ shown in Table 3 (left). Indicated by its coincidence matrix shown in Table 3 (center), the Robin Hood algorithm finally decides to swap the levels in the 4th factor between the first and last runs, in order to equalize pairwise coincidences. By such a single swap operation, we find that the PC-vector of the swapped design is majorized by the original PC-vector; see Figure 1 for the cumulative plots of sorted PC-vectors in the sense of (2.1).

In Figure 1 the benchmark by Theorem 1 is also plotted, which has two slopes (corresponding to $\theta$ and $\theta + 1$) rather than the dashed straight line, since $\bar{\beta} = 2.5714$ for $\mathscr{U}(8, 2^6)$. Note that the benchmark can be attained by any 6-factor sub-design of the 8-run Hadamard design. The Robin Hood swap algorithm makes the random design move toward the benchmark. Iterative swapping can make it move closer. However, the above deterministic procedure often gets into the local optimum. Advanced stochastic optimization methods are therefore called for. Based on a similar column-wise swap, Fang, Lu and Winker [7] used the threshold accepting heuristic for constructing uniform designs. Our group is currently developing a similar heuristic based on Robin Hood swaps, which is beyond the scope of this paper.

## APPENDIX

PROOF OF THEOREM 2. Let us write the sum notation $\sum_{u \subseteq \mathcal{S}, |u|=j}$, $\sum_{\mathbf{x} \in \mathcal{L}(q^j)}$ as $\sum_u$ and $\sum_{\mathbf{x}}$ in short, respectively. Let $\delta_{ik}^{(u)} = 1$ if the $u$-coordinate sub-tuples of $\mathbf{x}_i, \mathbf{x}_k$ take the same level combination and 0 otherwise. It can be verified that

$$\sum_u \delta_{ik}^{(u)} = \binom{\beta(\mathbf{x}_i, \mathbf{x}_j)}{j},$$

$$\sum_{i,k=1}^n \delta_{ik}^{(u)} = \sum_{\mathbf{x}} (N_{\mathbf{x}}^{(u)})^2,$$

and the Schur-combinatorial criterion (3.3) can be expressed as

$$\Psi_{\mathrm{C}}(\mathbf{X}; j) + \frac{n^2}{q^j}\binom{s}{j} = 2\sum_{r=1}^m \binom{\beta_r}{j} + n\binom{s}{j}$$

$$= \sum_{i,k=1}^n \binom{\beta(\mathbf{x}_i, \mathbf{x}_j)}{j}$$

$$= \sum_{i,k=1}^{n} \sum_{u} \delta_{ik}^{(u)} = \sum_{u} \sum_{\mathbf{x}} (N_{\mathbf{x}}^{(u)})^2.$$

TABLE 3
*Robin Hood swap on $\mathbf{X}(8, 2^6)$: at $j_* = 4$, swap the levels between rows $i_* = 1$ and $t_* = 8$; only 12 entries (boldfaced) are updated in the coincidence matrix*

| | | | | | | | | Original $\mathbf{M}$ | | | | | | | | Updated $\mathbf{M}$ | | | | | | |
|---|---|---|---|---|---|---|---|---|---|---|---|---|---|---|---|---|---|---|---|---|---|---|
| 0 | 0 | 0 | **1** | 1 | 0 | $i_*$ | | **5** | 3 | 1 | 1 | 3 | 5 | 0 | | **4** | **4** | **2** | **2** | **2** | **4** | 0 |
| 1 | 0 | 0 | 1 | 1 | 0 | $k$ | | | 2 | 2 | 0 | 4 | 4 | 1 | | | 2 | 2 | 0 | 4 | 4 | **2** |
| 0 | 0 | 1 | 0 | 1 | 1 | | | | | 2 | 4 | 2 | 2 | 3 | | | | 2 | 4 | 2 | 2 | **2** |
| 1 | 1 | 0 | 0 | 0 | 1 | | | | | | 4 | 2 | 2 | 5 | $\Rightarrow$ | | | | 4 | 2 | 2 | **4** |
| 0 | 1 | 1 | 0 | 0 | 1 | | | | | | | 2 | 2 | 5 | | | | | | 2 | 2 | **4** |
| 1 | 0 | 1 | 1 | 0 | 0 | | | | | | | | 2 | 3 | | | | | | | 2 | **4** |
| 0 | 1 | 0 | 1 | 1 | 0 | | | | | | | | | 1 | | | | | | | | **2** |
| 1 | 1 | 1 | **0** | 0 | 1 | $t_*$ | | | | | | | | | | | | | | | | |

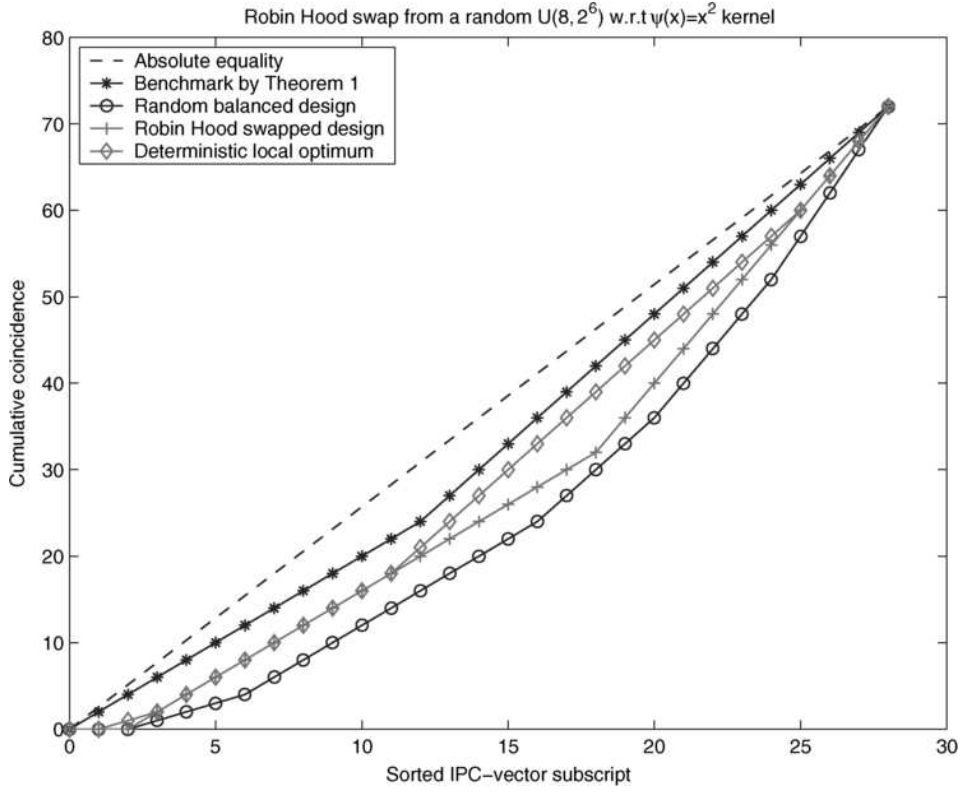

Fig. 1. *Robin Hood swap of a randomly generated* $\mathbf{X} \in \mathscr{U}(8, 2^6)$ *under kernel* $\psi(x) = x^2$.

By the method of variance decomposition, the right-hand side of (3.4) is given by

$$\sum_{u} \sum_{\mathbf{x}} \left(N_{\mathbf{x}}^{(u)} - \frac{n}{q^j}\right)^2 = \sum_{u} \sum_{\mathbf{x}} \left((N_{\mathbf{x}}^{(u)})^2 - \frac{2n}{q^j} N_{\mathbf{x}}^{(u)} + \frac{n^2}{q^{2j}}\right)$$

$$= \sum_{u} \sum_{\mathbf{x}} (N_{\mathbf{x}}^{(u)})^2 - \frac{n^2}{q^j} \binom{s}{j},$$

which equals the left-hand side $\Psi_C(\mathbf{X}; j)$. □

PROOF OF THEOREM 3. For the $j$-factor sub-design $\mathbf{X}_u(n, q^j)$ with $u$-coordinate factors, the word-length pattern (3.1) for $1 \le k \le j$ can be written as

$$A_k(\mathbf{X}_u) = \frac{1}{n^2} \sum_{wt(\mathbf{v})=k} |\chi_{\mathbf{v}}(\mathbf{X}_u)|^2,$$

where $\{\chi_{\mathbf{v}}, \mathbf{v} \in \mathcal{L}(q^j)\}$ are given orthonormal contrasts and $wt(\mathbf{v})$ is the number of nonzero elements of $\mathbf{v}$ [26]. Consider the deviation (3.5) of sub-design $\mathbf{X}_u$, $B_j^2(\mathbf{X}_u) = \frac{1}{q^j} \sum_{i,k=1}^n \binom{\beta_{ik}^{(u)}}{j} - \frac{n^2}{q^{2j}}$, where $\beta_{ik}^{(u)}$ measures the coincidence between $u$-coordinate sub-tuples of $\mathbf{x}_i, \mathbf{x}_k$ and it cannot exceed $j$. So $B_j^2(\mathbf{X}_u) = \frac{n}{q^j} E_0(\mathbf{X}_u) - \frac{n^2}{q^{2j}}$, where $E_0(\mathbf{X}_u)$ is defined in (3.2). Since it is true that

$$(3.6a) \qquad E_0(\mathbf{X}_u) = \frac{n}{q^s} \left(1 + \sum_{k=1}^{j} A_k(\mathbf{X}_u)\right)$$

(verified at the end), it follows that $B_j^2(\mathbf{X}_u) = \frac{n^2}{q^{2j}} \sum_{k=1}^{j} A_k(\mathbf{X}_u)$.

For $\mathbf{X}(n, q^s)$ itself, $B_j(\mathbf{X}) = \sqrt{\sum_{|u|=j} B_j^2(\mathbf{X}_u)}$. Via $\{\chi_{\mathbf{v}}, \mathbf{v} \in \mathcal{L}(q^j)\}$,

$$B_j^2(\mathbf{X}) = \frac{n^2}{q^{2j}} \sum_{|u|=j} \sum_{k=1}^{j} A_k(\mathbf{X}_u) = \frac{1}{q^{2j}} \sum_{|u|=j} \sum_{k=1}^{j} \sum_{wt(\mathbf{v})=k} |\chi_{\mathbf{v}}(\mathbf{X}_u)|^2.$$

Note that each contrast $\chi_{\mathbf{v}}(\mathbf{X}_u)$ of the sub-design $\mathbf{X}_u$ is also a contrast $\chi_{\mathbf{w}}(\mathbf{X})$ of $\mathbf{X}(n, q^s)$, where $\mathbf{w}$ coincides with $\mathbf{v}$ at $u$-coordinates and has null elements elsewhere. Denote by $\Omega_u$ the set of such $\mathbf{w}$'s. Then we have

$$B_j^2(\mathbf{X}) = \frac{1}{q^{2j}} \sum_{k=1}^{j} \sum_{|u|=j} \sum_{\mathbf{w} \in \Omega_u} |\chi_{\mathbf{w}}(\mathbf{X})|^2$$

$$= \frac{1}{q^{2j}} \sum_{k=1}^{j} \binom{s-j}{j-k} \sum_{wt(\mathbf{w})=k} |\chi_{\mathbf{w}}(\mathbf{X})|^2.$$

By writing $\frac{1}{n^2}\sum_{wt(\mathbf{w})=k}|\chi_\mathbf{w}(\mathbf{X})|^2$ back to $A_k(\mathbf{X})$, (3.6a) is proved. The benchmark of deviation pattern follows from Theorem 1 directly. Since the word-length patterns are linearly related through (3.6a) with positive pivoting coefficients, both patterns are equivalent under $\models$. We can therefore use $\widetilde{\boldsymbol{\beta}}$ that determines the benchmark of deviation pattern to derive the benchmark for the word-length pattern.

Let us now verify (3.6a) for $\mathbf{X}(n,q^s)$ through Krawtchouk polynomials. For an integer $l$ ($0 \leq l \leq s$) and any real number $y$, Krawtchouk polynomials $P_j(l;s,q)$ have the following property: $\sum_{j=0}^s P_j(l;s,q)y^j = [1+(q-1)y]^{s-l}(1-y)^l$ ([18], Section 5.7). By setting $y=1$, we have $\sum_{j=0}^s P_j(0;s,q) = q^s$ and $\sum_{j=0}^s P_j(l;s,q) = 0$, for $l = 1, \ldots, s$. By $P_0(l;s,q) = 1$ and $\sum_{l=0}^s E_0(\mathbf{X}) = n$,

$$\begin{aligned}
n\left(1 + \sum_{j=1}^s A_j(\mathbf{X})\right) &= n + \sum_{j=1}^s \sum_{l=0}^s E_l(\mathbf{X}) P_j(l;s,q) \\
&= n + \sum_{l=0}^s E_l(\mathbf{X}) \left(\sum_{j=0}^s P_j(l;s,q) - 1\right) \\
&= n + q^s E_0(\mathbf{X}) - \sum_{l=0}^s E_l(\mathbf{X}) = q^s E_0(\mathbf{X}).
\end{aligned}$$
□

PROOF OF THEOREM 5. In Definition 2, there is no risk in letting $\widehat{\mathbb{K}}_\varnothing = 1$ and $D_\varnothing(\xi; \widehat{\mathbb{K}}_\varnothing) = 0$. By the expansion of tensor products and coincidence measurements,

$$\begin{aligned}
D^2(\mathbf{X}; a \vee b) &= -\sum_{\varnothing \subseteq u \subseteq \mathcal{S}} \mu^{|u|} + \frac{1}{n^2} \sum_{\varnothing \subseteq u \subseteq \mathcal{S}} \sum_{i,k=1}^n \prod_{j \in u} [b + (a-b)\delta_{x_{ij},x_{kj}}] \\
&= -(1+\mu)^s + \frac{1}{n^2} \sum_{i,k=1}^n \prod_{j=1}^s [1 + b + (a-b)\delta_{x_{ij},x_{kj}}] \\
&= -(1+\mu)^s + \frac{1}{n^2} \sum_{i,k=1}^n (1+a)^{\beta(\mathbf{x}_i,\mathbf{x}_k)}(1+b)^{s-\beta(\mathbf{x}_i,\mathbf{x}_k)} \\
&= \frac{2}{n^2} \sum_{r=1}^m \left(\frac{1+a}{1+b}\right)^{\beta_r} + \frac{(1+a)^s}{n} - (1+\mu)^s.
\end{aligned}$$

By letting $\rho = (1+a)/(1+b)$, the identity (3.6g) follows. Provided that $b < a$ and $a < q-1$, and the exponential base $1 < \rho < \infty$, we can use Theorem 1 to get the lower bound for categorical discrepancy. □

**Acknowledgments.** We thank both referees for many helpful suggestions and appreciate the helpful discussions with Professor Rahul Mukerjee, Professor Fred Hickernell, Professor Min-Yu Xie, Dr. Min-Qian Liu and Dr. Yu Tang.

A. ZHANG
DEPARTMENT OF STATISTICS
UNIVERSITY OF MICHIGAN
ANN ARBOR, MICHIGAN 48109-1092
USA
E-MAIL: ajzhang@umich.edu

R. LI
DEPARTMENT OF STATISTICS
  AND THE METHODOLOGY CENTER
PENNSYLVANIA STATE UNIVERSITY
UNIVERSITY PARK, PENNSYLVANIA 16802-2111
USA
E-MAIL: rli@stat.psu.edu

K.-T. FANG
A. ZHANG
DEPARTMENT OF MATHEMATICS
HONG KONG BAPTIST UNIVERSITY
KOWLOON TONG
HONG KONG
E-MAIL: ktfang@math.hkbu.edu.hk

A. SUDJIANTO
RISK MANAGEMENT QUALITY
  AND PRODUCTIVITY
BANK OF AMERICA
CHARLOTTE, NORTH CAROLINA 28255-0001
USA
E-MAIL: Agus.Sudjianto@bankofamerica.com